\documentclass{amsart}
\usepackage{amsfonts,amssymb,amscd,amsmath,enumerate,verbatim,calc}

\newcommand{\TFAE}{The following conditions are equivalent:}
\newcommand{\CM}{Cohen-Macaulay}

\newcommand{\wrt}{with respect to}

\newcommand{\aF}{\mathbf{a} }

\newcommand{\n}{\mathfrak{n} }
\newcommand{\m}{\mathfrak{m} }
\newcommand{\C}{\mathcal{C}_{\bullet} }
\newcommand{\D}{\mathcal{D}_{\bullet} }
\newcommand{\q}{\mathfrak{q} }
\newcommand{\R}{\mathfrak{r} }

\newcommand{\rt}{\rightarrow}
\newcommand{\xar}{\longrightarrow}
\newcommand{\ov}{\overline}

\newcommand{\wt}{\widetilde }

\newcommand{\image}{\operatorname{image}}

\newcommand{\reg}{\operatorname{reg}}
\newcommand{\grade}{\operatorname{grade}}
\newcommand{\depth}{\operatorname{depth}}

\newcommand{\red}{\operatorname{red}}

\theoremstyle{plain}

\newtheorem{thm}{Theorem}

\newtheorem{theorem}{Theorem}[section]
\newtheorem{corollary}[theorem]{Corollary}
\newtheorem{lemma}[theorem]{Lemma}
\newtheorem{proposition}[theorem]{Proposition}
\newtheorem{question}[theorem]{Question}

\theoremstyle{definition}
\newtheorem{definition}[theorem]{Definition}
\newtheorem{s}[theorem]{}
\newtheorem{remark}[theorem]{Remark}
\newtheorem{example}[theorem]{Example}
\newtheorem{obs}[theorem]{Observation}

\theoremstyle{remark}

\numberwithin{equation}{theorem}

\begin{document}

\title[fiber cone]{The Hilbert coefficients of the fiber cone   and \\ the $a$-invariant of the associated graded ring}

\author{Clare D'Cruz}
 \address{Chennai Mathematical Institute, Plot H1, SIPCOT IT Park
Padur PO, Siruseri 603103, India}
\email{clare@cmi.ac.in}

\author{Tony~J.~Puthenpurakal}
\address{Department of Mathematics, IIT Bombay, Powai, Mumbai 400 076, India}
\email{tputhen@math.iitb.ac.in}
\thanks{The second author was partly supported by IIT Bombay seedgrant 03ir053}
\date{\today}
\begin{abstract}
Let $(A,\m)$ be a Noetherian local ring with infinite residue
 field and let $I$ be an ideal in $A$ and let $F(I) = 
\oplus_{n \geq 0}I^n/\m I^n$ be the fiber cone of $I$. 
We prove certain relations among the Hilbert coefficients $f_0(I),f_1(I), f_2(I)$ of $F(I)$
 when the $a$-invariant of the associated graded ring $G(I)$ is negative. 
\end{abstract}

 \maketitle
\section{introduction}
Let $(A,\m)$ be a Noetherian local ring with
 \emph{infinite residue field} $k = A/\m$. Let $I$ be an ideal in
 $A$. The \emph{fiber cone} of $I$ is the standard graded 
$k$-algebra
$F(I) = \bigoplus_{n\geq 0}I^n/\m I^{n}.$
Set $l(I) = \dim F(I)$,  
the \emph{analytic spread} of $I$.
 The Hilbert polynomial of $F(I)$ is denoted by $f_I(z)$.  Write
$f_I(z) = \sum_{i = 0}^{l-1}(-1)^i f_i(I)\binom{z+l-1-i}{l-1-i}$ where $l = l(I)$
We call $f_i(I)$  the $i^{th}$ \emph{fiber coefficient} of $I$. 

Most recent results in the study of fiber cone
 involve the depth of $G(I)   = \bigoplus_{n \geq 0}I^n/ I^{n+1} $,
the associated graded ring of $I$. 
When $I$ is $\m$-primary there has been some research relating $f_0(I)$ (the
\emph{ multiplicity} of $F(I)$) with various other
invariants of $I$ (see \cite[4.1]{JayV05}, \cite[4.3]{CVP032} and \cite[3.4]{Cor05}). 
In the case of $G(I)$ the relations among the Hilbert coefficients  $e_0(I),e_1(I),e_2(I)$ are well known
(see  \cite{VaSix}). However there is no result relating $f_0(I), f_1(I)$ and $f_2(I)$. 
 The reason for this is not difficult
to find: any standard $k$-algebra can be thought as a fiber cone
of its graded maximal ideal. So any result involving the 
relation between $f_i(I)$  \emph{would only hold in a restricted class of
 ideals}. 
Our  paper 
explores the relation between $\aF(I)$, \emph{the  $a$-invariant of $G(I)$},  and the
 Hilbert coefficients of $F(I)$. 
This is a new idea.

We first analyze when  
$l(I) = 2,3$ as it throws light on the general result.
 \begin{thm}\label{result1}
Let $(A,\m)$ be a Noetherian local ring with infinite residue field $k = A/\m$. Let $I$ be an ideal with $l(I) =2$. If $\aF(I) < 0$
 then 
\begin{center}
$f_1(I) \leq f_0(I) -1. $
\end{center}
 Furthermore,  equality holds if and only if $F(I^n)$ is \CM \ for all $n \gg 0$. If $\grade (I) = 2$ then equality holds.
\end{thm}
This result should be compared with a result due to Northcott \cite{Nor-e1}, which in our context states that
$f_1(\m) \geq f_0(\m) -1$ whenever $A$ is \CM. In \ref{2ex1}, we give an example of a two dimensional Noetherian local ring $(A,\m)$
with $\depth A = 1$ but $f_1(I) < f_0(I) -1$.

To analyze the case when equality holds  in Theorem \ref{result1} we resolve $F(I^n)$ as a $F(J^{[n]})= k[X^{n}_{1}, X^{n}_{2}]$-module 
and write
it   as:
\[
0 \xar K_n \xar \bigoplus_{i = 1}^{\beta_{1}^{[n]}}F(J^{[n]})( - 1 - \alpha_{i}^{[n]} )  \xar F(J^{[n]})^{\beta_{0}^{[n]}} \xar F(I^n) \xar 0
\]
Here $\alpha_{i}^{[n]} \geq 0$ for all $i$.
As $ \depth F(I^n) \geq 1$ for all $n \gg 0$ we get
$K_n = 0$ for all $n \gg 0$.
We show in Theorem \ref{sub} that if $\aF(I) < 0$ then \emph{for all} $n \gg 0$,
\begin{align*}
f_1(I) - f_0(I) + 1 &= - \sum_{i = 1}^{\beta_{1}^{[n]}}\alpha_{i}^{[n]}   \quad \text{and} \\
\beta_{1}^{[n]} &= 0 \ \ \text{if and only if } \alpha_{i}^{[n]} = 0 \ \text{ for all $i$ }.
\end{align*}

Our second result Theorem \ref{sectheorem} has a 
noteworthy consequence when $G(I)$ is \CM.

\begin{thm}\label{seccor}
Let $(A,\m)$ be \CM \ local ring  of dimension $d = 3$. Let  $I$ be an $\m$-primary
ideal with $G(I)$ \CM \ and  $\red(I) = 2$. Then
\begin{center}
$f_2(I) \geq f_1(I) - f_0(I) + 1$.
\end{center}
\end{thm}

We extend our results to higher analytic spread using Rees-superficial 
sequences (see section $6$  for details), under some mild assumptions on  $grade(I)$.
We state some of our noteworthy results. The first one (see \ref{mth1}) states that
if $l(I) \geq 2 $, $\grade(I) \geq l(I) -2$  and reduction number
of $I \leq 1$ then $f_1(I) \leq f_0(I) - 1$ with equality if $\grade(I) = l(I)$. 
An immediate consequence
(see \ref{NarC})  is that if $(A,\m)$ is \CM \ with $\dim A \geq 2$, $I$ an $\m$-primary ideal  and  the second Hilbert-Samuel coefficient $e_2(I) = 0$ 
then $f_1(I) = f_0(I) - 1$ (see 2.7 for definition of $e_2(I)$).

Finally, we show that if $A$ is \CM \ ring of dimension at least three
and if  $I$  an $\m$-primary ideal of reduction number two whose associated graded ring is  \CM,  then 
$f_2(I) \geq f_1(I)-f_0(I) +1$ (see \ref{mth2}).

 Here is an overview of the contents of the paper. 
In section $1$
we introduce some notations and  preliminary facts needed. 
 In section $2$ we introduce two complexes which will be used in the subsequent 
sections.
In section $3$ we prove the main result for $l=2$ 
(Theorem \ref{result1}).
 In
 section $4$
we prove our second Theorem and as a consequence obtain Theorem \ref{seccor}.
 In section $5$ we obtain  results on the coefficients of the fiber cone 
for   any analytic spread. In the appendix( =section 6)  we recall some basic facts regarding minimal reductions and
filter-regular elements and prove an elementary result; which is useful in section 3. 

\textit{Acknowledgments}
The authors thanks the referee for many pertinent
 comments. The author also thanks Fahed Zulfeqarr and A. V. Jayanthan for help in examples. 
 
\section{preliminaries}
 From now on
$(A,\m)$ is a Noetherian local ring of dimension $d$, with infinite residue field. 
All modules are assumed to be finitely
generated.
 For a  finitely generated module  $M$, we  denote its length by
$\ell(M)$.

\s Let $J = (x_1,\ldots,x_l)$ be  a minimal reduction of 
$I$.  We denote by  $\red_J(I) := \min \{ n | JI^n= I^{n+1} \} $ 
the \textit{reduction number} of $I$ \wrt \ $J$. Let
 $$\red(I) = \min \{ \red_J(I) \mid J
\text{\ is  a reduction of \ } I \}$$
 be the \textit{reduction number } of $I$.

 \s As a reference for local cohomology we use  \cite{BSh} (see especially Chapter 18 for relations
between local cohomology and reductions). 
We take local cohomology of $G(I)$ \wrt \ $G(I)_+ =  \bigoplus_{n\geq 1}I^n/I^{n+1}$. Set $G = G(I)$ and $G_+ = G(I)_+$.  For each $i \geq 0$ 
the local cohomology modules
$H^i_{G_+}(G)$ are  graded   and furthermore
 $H^i_{G_+}(G)_n = 0$ for all $n \gg 0$.
For each $i \geq 0$ set
\[
a_i(I) = \max \{ n \mid H^i_{G_+}(G)_n \neq 0 \}.
\]

\s  Set $l = l(I)$. Then
$H^l_{G_+}(G) \neq 0$ and $H^i_{G_+}(G) = 0 $ for all $i > l$ (see \cite[2.3]{Hoa}).
We call $ \aF(I) = a_l(I)$ to be the \emph{$a$-invariant} of $G(I)$.
The \emph{(Castelnuovo-Mumford) regularity} of  $G(I)$ is
\[
\reg(G(I)) =  \max \{ a_i(G) + i \mid 0 \leq i \leq l \}.
\]
The \emph{regularity of $G(I)$ at and 
above level $r$}, denoted by $\reg^r(G(I))$, is
\[
\reg^r(G(I)) =  \max \{ a_i(G) + i \mid r \leq i \leq l \}.
\]

\s\label{notation}
Let $x \in I\setminus I^2$ be a  $I$-superficial element of $I$. For all
 $r \geq 1$ and  $s \geq 0$ it is  easy to show
\[
\reg^{r}(G(I)) \leq s \implies \reg^{r}(G(I/(x)) \leq s.
\]

\s\label{asymp} We will use the following beautiful result due to Hoa(\cite[2.6]{Hoa}):

\noindent\emph{ There exists non-negative integers $n_0, \R(I)$, such that for all $n \geq n_0$ and every minimal reduction $J$ of $I^n$  we have
 $\red_J(I^n)  = \R(I)$. Furthermore}
\begin{equation*}
\R(I) =
\begin{cases}
l(I) - 1 & \text{if } \ \aF(I) < 0, \\
l(I) & \text{if } \ \aF(I) \geq 0.
\end{cases}
\end{equation*}

\s For definition and  basic properties of superficial sequences see \cite[p.\ 86-87]{Pu1}
.
\s  If $I$ is $\m$-primary then let  $p_I(z)$ be the \emph{Hilbert-Samuel polynomial} of $A$ \wrt \  $I$ 
(so  $\ell (A/I^{n+1}) = p_I(n) $, $\forall n \gg 0$).
Write
$
p_I(z) = \sum_{i = 0}^{d}(-1)^i e_i(I)\binom{z+d -i}{d-i}.
$
For $i \geq 0$ we call $e_i(I)$  the $i^{th}$-\emph{Hilbert coefficient} of $I$.

\s The Hilbert series of $F(I)$, $G(I)$ is denoted by $H(F(I),z)$, $H(G(I),z)$ respectively,  i.e.
\[
H(F(I),z) = \sum_{n \geq 0} \ell\left( \frac{I^n}{\m I^{n}} \right)z^n\quad \text{and} \quad H(G(I),z) = \sum_{n \geq 0} \ell\left( \frac{I^n}{I^{n+1}}  \right)z^n
\]

\s  If $x \in I^j \setminus \m I^j$ then we denote by $x^\circ$ its image in $F(I)_j = I^j/\m I^j$.

\section{Two complexes}
Our results are based on analyzing two complexes which 
we describe in this section. Both the complexes are defined 
\emph{using the maps in the Koszul complex.} 
Throughout $I$ is an ideal in $(A, \m)$ with a minimal reduction $J = (x_1,\ldots,x_l)$ where $l =l(I)$. Note that
$x_1,\ldots,x_l$ are analytically independent (\cite{NoRees}, cf. \cite[4.6.9]{BH}). Using analyticity to 
detect exactness (at some stage) of a complex

\s The first complex is when $l = 2$
\begin{equation}
\label{first}
\mathcal{C}_\bullet(I): \quad 0\rt \frac{A}{\m} \xrightarrow{\alpha_2} \left(\frac{I}{\m I}\right)^2\xrightarrow{\alpha_1}  \frac{IJ}{\m IJ} \rt 0, \quad \text{where} 
\end{equation}
\[
     \alpha_2(a + \m) 
= \begin{pmatrix}-x_2a + \m I \\ x_1a + \m I\end{pmatrix} \ \  \& \ \ 
\alpha_1 \begin{pmatrix}a + \m I \\ b + \m I\end{pmatrix} = x_1a + x_2b + \m IJ.
\]
\begin{obs}\label{firstcomments}
\

\begin{enumerate}[\rm(i)]
\item
$H_2(\C(I)) = 0$ since  $x_1,x_2$ are analytically independent.
\item
If $x_1,x_2$ is a regular sequence in $A$
then clearly
$H_1(\C(I))  = 0$.  
\item
Clearly
  $\alpha_1$ is surjective. So $H_0(\C(I)) = 0$.
\end{enumerate}
\end{obs}

\s The second complex, $\mathcal{D}_\bullet(I)$, is when $l = 3$.
\begin{equation}
\label{second}
     0
\rt                          \frac{A}{\m} 
\xrightarrow{\alpha_3} \left(\frac{I}{\m I}\right)^3
\xrightarrow{\alpha_2} \left(\frac{I^2}{\m I^2}\right)^3  
\xrightarrow{\alpha_1} \frac{I^2J}{\m I^2J} 
\rt                      0, \quad \text{where}
\end{equation}
\[
  \alpha_3(a + \m) 
= \begin{pmatrix}
    x_3a + \m I \\
   -x_2a + \m I \\ 
    x_1a + \m I
  \end{pmatrix}, \quad 
  \alpha_1
  \begin{pmatrix}
   a + \m I^2 \\ 
   b + \m I^2 \\ 
   c + \m I^2\end{pmatrix} 
= x_1a + x_2b + x_3c + \m I^2J,
\]
\[
\text{and}\quad \alpha_2\begin{pmatrix}a + \m I \\ b + \m I \\ c + \m I\end{pmatrix} = 
\begin{pmatrix}-x_2a -x_3b + \m I^2 \\ x_1a -x_3c + \m I^2 \\  x_1b + x_2c + \m I^2\end{pmatrix}.
\]
\begin{obs}\label{secondcomments}
\

\begin{enumerate}[\rm (i)]
\item
$H_3(\D(I)) = 0$ since  $x_1,x_2, x_3$ are analytically independent.
\item
In Lemma \ref{2c1} we show that if $x_1,x_2,x_3$ is a regular sequence and if
$I^2\cap J = JI$ then $\image(\alpha_2) = \ker (\alpha_1)$, so
$H_1(\D(I))  = 0$.  
\item
The assumption $I^2\cap J = JI$ holds when the following hold
\begin{enumerate}
\item
 $I$ is integrally closed (see\cite[p.\ 317]{Hun} and \cite[Theorem 1]{Itoh88}). 
\item
 The initial forms $x_{1}^{*}, x_{2}^{*}, x_{3}^{*}$ in $G(I)_1$ is a
regular sequence (see \cite[2.3]{VV}).
\end{enumerate}
\item
Clearly
$\alpha_1$ is surjective and so $H_0(\D(I)) = 0$.
\end{enumerate}
\end{obs}

\begin{lemma}\label{2c1}(with hypothesis as in \ref{second})
If $x_1,x_2,x_3$ is a regular sequence and if $I^2\cap J = JI$ then $\image(\alpha_2) = \ker (\alpha_1)$.
\end{lemma}
\begin{proof}
Let $\mathcal{K}(\mathbf{x})_\bullet$ be the Koszul complex on 
$x_1,x_2,x_3$. It is acyclic since  $x_1,x_2,x_3$ is a regular sequence.
Suppose 
%$ \xi = [\ov{a},\ov{b},\ov{c}]^t \in \ker\alpha_1$.
$ \xi 
= \begin{pmatrix}
   \ov{a} \\
   \ov{b} \\
   \ov{c}
\end{pmatrix}
 \in \ker\alpha_1
$.
 Then
\[
ax_1 + bx_2 + cx_3 = a'x_1 + b'x_2 + c'x_3 \quad \text{where} \ a',b',c' \in \m I^2.
\]
As  $\mathcal{K}(\mathbf{x})_\bullet$ is acyclic;  there exists $f,g,h \in A $ such that

\begin{equation}\label{lemmeq}
\begin{pmatrix}a -a' \\ b-b' \\ c-c' \end{pmatrix} = \begin{pmatrix}-x_2f -x_3g \\ x_1f -x_3h \\  x_1g +
 x_2h \end{pmatrix}.
\end{equation}
We  show that $f,g,h$ are in $I$. Using (\ref{lemmeq}) we get
\[
-x_2f -x_3g \in I^2\cap J = JI.
\]
So $-x_2f - x_3g = x_1p + x_2q + x_3r$ where $p,q,r \in I$. Again using the fact that 
 $\mathcal{K}(\mathbf{x})_\bullet$ is acyclic we get that there exists $u,v,w \in A$ such that
\[
\begin{pmatrix}p  \\ q+ f \\ r+ g \end{pmatrix} = \begin{pmatrix}-x_2u -x_3v \\ x_1u -
x_3w \\  x_1v +
 x_2w \end{pmatrix}.
\]
So $f,g \in I$. Similarly by using the
 second row in (\ref{lemmeq})  we get
that $f,h \in I$.  
Set 
%$\eta = [\ov{f},\ov{g},\ov{h}]^t \in (I/\m I)^3$. 
$  \eta 
=  \begin{pmatrix}
   \ov{f} \\
   \ov{g} \\
   \ov{h}
\end{pmatrix}
 \in (I/\m I)^3$. 
Notice
$\alpha_2(\eta) = \xi $ since $a', b', c' \in \m I^2$.
 Thus
$\ker \alpha_1 \subseteq \image \alpha_2$.
\end{proof}

\s Let us recall the following well-known fact of complexes:
Let 
\[
\mathcal{X}_{\bullet}: \quad 0\rt X_n \rt X_{n-1}\rt \ldots \rt X_0 \rt 0,
\]
be a complex of $A$-modules with $\ell(X_i)$ finite for all $i$. Then
\begin{equation}\label{hom}
\sum_{i = 0}^{n}(-1)^i\ell(X_i) = \sum_{i = 0}^{n}(-1)^i\ell( H_i\big(\mathcal{X}_{\bullet}\big)).
\end{equation}

\section{Proof of Theorem \ref{result1}}
 In this section we prove Theorem \ref{result1}. The setup below is 
used throughout. The hypothesis $\aF(I)$ is crucial for our results. See
 examples 6.1, 6.2,6.3 for some illustrations of Theorem \ref{result1}.

\s \textbf{Setup:} Let $J = (x_1,x_2)$ be a minimal reduction of $I$.
Notice $ J^{[n]}= (x_{1}^{n},x_{2}^{n})$ is a minimal reduction of $I^n$.
 If $\grade(I) = 2$ then we can take $x_1,x_2$
to be a regular sequence (and so $x_{1}^{n},x_{2}^{n}$ is also a regular sequence).
For $i = 1, 2$ set $X_i = x_{i}^{\circ}$, the image of $x_i$ in $I/\m I$. So
$F(J^n) = k[X^{n}_{1}, X^{n}_{2}]$ for all $n \geq 1$.
 
We first prove the inequality stated in Theorem  \ref{result1}.
\begin{theorem}\label{1}
Let $(A,\m)$ be a  local ring and let $I$ be an ideal with $l(I) =2$. If $\aF(I) < 0$
 then 
$f_1(I) \leq f_0(I) -1.$
Furthermore if $\grade(I) = 2$ then  equality holds.
\end{theorem}
\begin{proof}
We consider the complex (\ref{first}) for the ideal $I^n$ for each $n \geq 0$. 
Set $\C [n] = \C (I^n)$ for $n \geq 1$. By Remark \ref{firstcomments} we have
$H_i(\C [n]) = 0$ for $i = 0,2$.
Using (\ref{hom}) for the complex $\C[n]$ for each $n$ we get an equation
\begin{equation}\label{r1}
1 - 2\ell(I^n/\m I^n) + \ell(I^nJ^{[n]}/\m I^nJ^{[n]}) = -\ell(H_1(\C [n])).
\end{equation}

 Since $\aF (I) < 0$,  by \ref{asymp} we have $\red_{J^{[n]}}(I^n) =1$ for all $n \gg 0$. So
$I^nJ^{[n]} = I^{2n}$ for all $n \gg 0$. Also for all $n \gg 0$ we have
$f_I(n) = \ell(I^n/\m I^n)$.
Setting these in (\ref{r1})
we get for all $n \gg 0$,
\[
1 - 2f_I(n) + f_I(2n) = -\ell(H_1(\C [n])).
\]
Write $f_I(n) = f_0(n+1) -f_1$. Therefore
\[
1 - 2\{f_0(n+1) -f_1 \} + f_0(2n+1) -f_1 = -\ell(H_1(\C [n])).
\]
\[
\text{Thus} \quad 1 - f_0 + f_1  =  -\ell(H_1(\C [n])).
\]
Hence
\[
\text \quad 1 - f_0 + f_1  \leq 0.
\]
By our  assumption on $J$,  $H_1(\C [n]) = 0$ for each $n$ if $\grade(I) = 2$.
Hence equality holds in the above equation.
This proves the result. 
\end{proof}

The following example shows that if $\aF(I) < 0$ but $\grade(I) \neq 2$ then $f_1(I) < f_0(I) -1 $
 is possible. 
\begin{example} \label{2ex1}
Set $A = k[[X_1,X_2, X_3]]/(X_{1}^{2}, X_1X_2) = k[[x_1,x_2,x_3]]$. Set $I = \m = (x_1,x_2,x_3)$.
and $J = (x_2,x_3)$. Then $J$ is a reduction of $\m$ and  $\m^2 = J\m$. By \cite[3.2]{TrungProc} we get $\aF(\m)< 0$.

It can be checked that $\grade (\m) = 1$ and $x_3$ is a non-zero divisor. 
Using COCOA \cite{cocoa}
 it can verified that the Hilbert-series of $F(\m) = G(\m)$ is
\[
\frac{1+z-z^2}{(1-z)^2}.
\]
So $f_1(\m) = -1$ but $f_0(\m) = 1$.
\end{example}

The next example shows that Theorem \ref{1}  need not hold  when $  \aF(I)> 0$.
\begin{example}\label{2ex2}
Let
$(A,\m)$ be a two dimensional \CM \ local ring  with $\red(\m) = 2$. Then will have
$G(\m) = F(\m)$ is  \CM, (\cite[2.1]{Sa2}) and   it's  Hilbert-series is
\[
\frac{1+z+ cz^2}{(1-z)^2} \quad \text{where} \ c > 0.
\]
So $f_1(\m)-f_0(\m) + 1 = c > 0$.  
\end{example}

Next we analyze the case when $f_1(I) = f_0(I) -1$. 
\s \label{veronese}Observe that $F(I^n) = F(I)^{<n>}$, the 
 \emph{$n^{th}$ Veronese} subring of $F(I)$. In particular $l(I^n) = l(I)$ for each $n \geq 1$. Local cohomology commutes with
the Veronese functor, cf. \cite[2.5]{HMc}. In \cite[2.8]{Pu3} it is proved that
$\depth F(I^n)$ is constant for all $n \gg 0$. We prove

\begin{theorem}\label{veronese-th}
 Let $(A,\m)$ be a local ring  with infinite residue field and let $I$ be an ideal with $s = l(I) > 0$. Then
\begin{enumerate}[\rm (1)]
 \item 
$\depth F(I^n) > 0$ for all $n \gg 0$.
\item 
There exists a minimal reduction $J = (x_1,\ldots,x_s)$ of $I$ such that $(x_1^n)^\circ$ is $F(I^n)$ regular for 
all $n \gg 0$. 
\end{enumerate}
\end{theorem}
\begin{proof}
 (1) Set $E = H^{0}_{F(I)_+}(F(I))$. Clearly $\ell(E) < \infty$. Say $E = \bigoplus_{i = 0}^{r}E_i$.  Notice $E$ is an
ideal of $F(I)$ with finite length. If
$E_0 \neq 0$ then $1_{F(I)} \in E$.  So $E = F(I)$ will have finite length, a contradiction since $\dim F(I) = l(I) \geq 1$. Therefore $E_0 = 0$ and as a consequence we have
\[
 H^{0}_{F(I^n)_+} F(I^n) =  \left( H^{0}_{F(I)_+}(F(I) \right)^{<n>} = 0  \quad \text{for all} \ \ n > r.
\]
Thus $\depth F(I^n) > 0$ for all $n \gg 0$.

(2) By \cite[3.8]{TrungProc}  we get that 
 there exists a minimal reduction $J = (x_1,\ldots, x_s)$ of $I$ such that
$x_1^\circ, \ldots, x_s^\circ \in F(I)_1$ is a $F(I)$-filter regular sequence. Set $x = x_1$. Since $x^\circ$ is $F(I)$-filter regular, by \ref{filt-ideal-cor} we get that $(x^n)^\circ$ is $F(I^n)$-filter regular for each $n \geq 1$.
By (1)  $\depth F(I^n) > 0$ for all $n \gg 0$. So by \ref{filtISreg} we get that $(x^n)^\circ$ is $F(I^n)$-regular for all $n \gg 0$.
\end{proof}

\s \label{rossi}
As $l(I) = 2$, by computing the Hilbert polynomial of $F(I)$ and $F(I^n)$  we obtain 
\[
f_1(I^n) - f_0(I^n) + 1 = f_1(I) - f_0(I) + 1 \quad \text{for all } n \geq 1.
\]

\s \label{temp} Let $J = (x_1,x_2)$ be a minimal reduction of $I$ as constructed in Theorem \ref{veronese-th}. 
In particular $(x^{n}_{1})^\circ$ is $F(I^n)$-regular for all $n \gg 0$. Set $X^{n}_{j} = (x^{n}_{1})^\circ$ for $j = 1,2$.
 We resolve $F(I^n)$ as a $F(J^{[n]}) =  k[X^{n}_{1}, X^{n}_{2}]$ module and write
it   as:
\[
0 \xar K_n \xar \bigoplus_{i = 1}^{\beta_{1}^{[n]}}F(J^{[n]})( - 1 - \alpha_{i}^{[n]} )  \xar F(J^{[n]})^{\beta_{0}^{[n]}} \xar F(I^n) \xar 0
\]
Here $\alpha_{i}^{[n]} \geq 0$ for all $i$.  As $ \depth F(I^n) \geq 1$ for all $n \gg 0$ we get
$K_n = 0$ for all $n \gg 0$.
We prove

\begin{theorem}\label{sub}
(with assumptions as in \ref{temp}) If $\aF(I) < 0$ then
 \emph{for all} $n \gg 0$,
\begin{align*}
f_1(I) - f_0(I) + 1 &= - \sum_{i = 1}^{\beta_{1}^{[n]}}\alpha_{i}^{[n]}   \quad \text{and} \\
\beta_{1}^{[n]} &= 0 \ \ \text{if and only if } \alpha_{i}^{[n]} = 0 \ \text{ for all $i$ }.
\end{align*}
\end{theorem}
For the proof of this theorem we need the following:
\begin{lemma}\label{euclid}
Let $R = k[X]$ and let $S = \oplus_{n \geq 0}S_n$ be a standard $k$-algebra of
dimension $1$ and multiplicity $p+1$. Assume $S = Ru_1 + \ldots +Ru_m$ where  degree $u_i \leq 1$ and 
$u_1 = 1_{S}$. Then  $S$ has the following resolution over $R$
\[
0 \xar \bigoplus_{i = 1}^{q}R( -1  - \alpha_i  )\xar R \oplus R(-1)^{p+ q} \xar S \xar 0. 
\]
Also $S$ is free if and only if all $\alpha_i = 0$. 
\end{lemma}
\begin{proof}
By the hypothesis on $S$ and as $R$ is an Euclidean domain we get
\[
S = R \bigoplus R(-1)^{p} \bigoplus \left(  \bigoplus_{i = 1}^{q}\frac{R}{(X^{\alpha_i}) } (-1)    \right), \quad  \text{where} \  \alpha_i \geq 0.
\]
The result follows.
\end{proof}

\begin{proof}[Proof of Theorem \ref{sub}]
We choose $n_0$ such that for all $n \geq n_0$ we have $\depth F(I^n) \geq 1$ and
$\red_{J^{[n]}}(I^n) = 1$.  Fix $n \geq n_0$. Set $\alpha_i = \alpha_{i}^{[n]} $.
Set $F(J^{[n]}) = k[X^{n}_{1}, X^{n}_{2}]$. Since $\red_{J^{[n]}}(I^n) = 1$ it follows 
that $F(I^n)$ is generated as a $F(J^{[n]})$ in degrees $\leq 1$.

Notice that by construction, $X_1^{n}$ is a non-zero divisor on $F(I^n)$ (see \ref{temp}).
Set $R = F(J^{[n]})/(X^{n}_{1}) =  k[X^{n}_{2}]$ and $S = F(I^n)/X^{n}_{1}F(I^n)$.
Note that  $S$ is generated as a $R$ module in degrees $\leq 1$. By Lemma \ref{euclid} 
the resolution of $S$ as a $R$-module is
\[
0 \xar \bigoplus_{i=1}^{q}R( -1  - \alpha_i)   \xar R \oplus R(-1)^{p+ q} \xar S \xar 0.
\]
Since $X_{1}^{n}$ is a non-zero divisor on $F(I^n)$ and $F(J^{[n]})$, we get
that the  resolution of $F(I^n)$ as  $F(J^{[n]})$-module is
\[
0 \xar \bigoplus_{i=1}^{q} F(J^{[n]})( -1 - \alpha_i) \xar F(J^{[n]}) \oplus F(J^{[n]})(-1)^{p+ q} \xar F(I^n) \xar 0.
\]
Thus $q = \beta_{1}^{[n]} $. Set $\phi(z) = \sum_{i = 1}^{q} z^{\alpha_i + 1}$. 
Therefore  the Hilbert series of $ F(I^n)$ is 
\[
  \frac{ 1 + (p+q)z -\phi(z)}{(1-z)^2}.
\]
So  $f_0(I^n) = 1 + p$. Notice
\[
f_1(I^n) = p + q - \sum_{i = 1}^{q}(\alpha_i + 1)  = p - \sum_{i = 1}^{q}{\alpha_i}.
\]
Using \ref{rossi} we obtain
\[
\sum_{i = 1}^{q}{\alpha_i} = f_0(I^n) - f_1(I^n) - 1 = f_0(I) - f_1(I) - 1.
\]
Also by Lemma \ref{euclid}, $q = 0$ if and only if all $\alpha_i = 0$.
\end{proof}
In view of this result we are tempted to ask
\begin{question}
(with notation as in \ref{temp})
 Let $(A,\m)$ be a local ring  and let $I$ be a proper ideal with $l(I) = 2$. Is  $ \sum_{i = 1}^{\beta_{1}^{[n]}}{\alpha_i}^{[n]}$
 constant for all $n \gg 0$?
\end{question}

Finally we prove
\begin{proof}[Proof of Theorem \ref{result1}]
This follows from Theorems \ref{1}, \ref{sub}.
\end{proof}

We give an example which shows that in the case $l(I) = \grade(I) = 2$ and $\aF(I) < 0$ it is possible for $F(I)$  to be not \CM \ even though
$f_1(I) = f_0(I) -1$ and 
 $F(I^n)$
is \CM  \ for all $n \gg 0$. The example below was constructed by Marley \cite[4.1]{Mar} in his study of the associated graded ring $G(I)$.
\begin{example}
 Let $A = k[X,Y]_{(X,Y)}$ and let $I = (X^7, X^6Y, XY^6, Y^7) $. Using COCOA one verifies that 
$e_2(I) = 0$ and that the Hilbert series of the fiber cone $F(I)$ is
\[
 \frac{1 + 2z + 2z^2  + 2z^3 + 2z^4 + 2z^5 - 4z^6 }{(1-z)^2}
\]
From the Hilbert series it's clear that $f_0(I) = 7$,
$f_1(I) = 6$ but  $F(I)$ is not \CM.

\end{example}

\begin{remark}
It is also possible to prove  Theorem \ref{result1} directly from Theorem \ref{sub},
Observation \ref{rossi} 
and by using a result of Kishor Shah \cite[Theorem 1]{shah91}.
We kept Theorem \ref{1} in this section since  it gives 
the inequality $f_1(I) \leq f_0(I) - 1$ very easily and more importantly
gives us a natural way of trying to relate $f_i(I) $ for $i = 0,1,2$ which is new and important.  This
is done in our next section.
\end{remark}
\section{Results when analytic spread is three} 
In this section
we assume that $l(I) = 3$. If $J = (x_1,x_2,x_3)$ is a reduction of $I$, we also assume that
$x_1,x_2,x_3$ is a regular sequence.  The goal of this section is to prove
$f_2(I) \geq f_1(I) - f_0(I) + 1$ under suitable conditions on $I$.

\s\label{scrap} We consider the complex (\ref{second}) for the ideal $I^n$ for each $n \geq 0$. 
Set $\D [n] = \D (I^n)$ for $n \geq 1$. By Remark \ref{firstcomments} we have
$H_i(\D [n]) = 0$ for $i = 0,3$. 

To analyze the case when $H_1(\D [n]) $ is zero we make the following:
\begin{definition}
Let $I$ be an ideal and let $J = (x_1,x_2,x_3)$ be a minimal reduction of $I$. We say  the pair $(I,J)$ satisfy 
$V_{2}^{\infty}$ if 
\begin{equation*}
I^{2n}\cap J^{[n]} = J^{[n]}I^n \quad \text{for all} \ n\gg 0. 
\end{equation*}
This condition has been studied by Elias in \cite{Elasymp}. 
\end{definition}

\s \label{scrap2} By Lemma
\ref{2c1},  $H_1(\D [n]) = 0$ for all $n \gg 0$
\emph{ if} the pair $(I,J)$ satisfy $V_{2}^{\infty}$.

\s \label{hyp-satisfy}  When $\grade(I) = l(I)$ then using \ref{secondcomments}(iii) the hypothesis $V_{2}^{\infty}$ holds  when
\begin{enumerate}
\item
$I$ is asymptotically normal i.e., $I^n$ is integrally closed for all $n \gg 0$.
\item
The initial forms $x_{1}^{*},\ldots, x_{l(I)}^{*}$ in $G(I)_1$ form a
regular sequence. 
\end{enumerate}

We now state our second main theorem.

\begin{theorem}\label{sectheorem}
Let $(A,\m)$ be local and let $I$ be an ideal in $A$  with $l(I) = \grade(I) =3$. Let $J = (x_1,x_2,x_3)$ be a minimal reduction of $I$ and assume the pair $(I,J)$ satisfy $V_{2}^{\infty}$. If $\aF(I) < 0$ then
$f_2(I) \geq f_1(I) -f_0(I) + 1$. 
\end{theorem}
\begin{remark}
 The hypothesis of the theorem above is quite stringent. However they are necessary (see Examples \ref{ex1}, \ref{ex2}, \ref{ex3}). Also note that if $F(I)$ is \CM \ then inequality holds (infact it holds for any
\CM \ $k$-algebra).
\end{remark}
We give two examples where the condition of Theorem \ref{sectheorem} holds. 
\begin{example}\label{pos1}
 Let $(A,\m)$ be a three dimensional \CM  \ local ring with $I$ an $\m$-primary ideal with reduction number two and
$G(I)$ \CM.  Then the hypothesis of Theorem \ref{sectheorem} holds by \ref{hyp-satisfy}(2). This is shown in Theorem \ref{as=2}. Using Example 6.1 from \cite{jpv} we can construct an interesting  example of this kind   as follows:

Let $T = k[[t^6, t^{11}, t^{15}, t^{31}]]$, $K = (t^6, t^{11},
t^{31})$ and $L = (t^6)$. Then, it can easily be verified that 
 $K^3 = LK^2$. Since $K^2 \cap L = LK$, 
$G(K)$ is Cohen-Macaulay by a result of Valabrega and Valla \cite[2.3]{VV}. 
It can also be seen that $t^{37} \in \m
K^2$, but $t^{37} \notin \m LK$. Therefore $F(K)$ is
not Cohen-Macaulay by a criterion due to Cortadellas  and Zarzuela \cite[3.2]{CZ}.
One can verify that the Hilbert-Series of $F(K)$ is $(1 + 2z)/(1-z)$.

Let $R = k[[X,Y,Z,W]]$ and 
$\q = (y^2z-xw,x^2z^2-yw,x^3z-y^3,x^3yw-z^4,z^5-y^4w,xyz^3-w^2,y^5-wx^4,x^2y^3-z^3,x^5-z^2,x^4y^2-zw)$.  Set
$B = R/\q= k[[x,y,z,w]]$. Using COCOA, one can verify that $B \cong T$. Under this isomorphism
the ideal $(x,y,w)$  maps to $K$ and $(x)$ goes to $L$

Set $A = B[[U,V]]$. Clearly $A$ is \CM \ of dimension $3$. Set $I = (x,y,w, U, V)$. and $J = (x, U, V)$.
Clearly $J$ is a minimal reduction of $I$ and $I^3 = JI^2$. Furthermore
$G(I) \cong G(K)[U, V]$ and $F(I) \cong F(K)[U, V]$. So $G(I)$ is \CM \ while $F(I)$ has depth $2$. 
\end{example}
Be giving the next example we make the following:
\begin{remark}
If $A$ satisfies the condition of the theorem then
as   $\aF(I) < 0$ we get by \ref{asymp} that $\red_{J^{[n]}}(I^n) = 2$ for all $n \gg 0$.  So by using the
Valabrega-Valla criterion \cite[2.3]{VV} it follows that
 $\grade G(I^n) \geq 3$ for all 
$n \gg 0$.
\end{remark}
\begin{example}\label{pos2}
  Let $A = k[[X,Y,Z]]$ and let $\m = (X,Y,Z)$.  Let $I$ is an $\m$-primary ideal with $I^r = \m^s$ for some $s> r$ and 
$G(I)$ not \CM \ ( for a specific example see \cite[4.3]{CPR} also see \cite[3.8]{CPR}).  If
$J$ is a reduction of $I$ then the pair $(I,J)$ satisfy $V_{2}^{\infty}$, by \ref{hyp-satisfy}(1).
  Also by Hoa's result it follows that $a(I) < 0$. Any such example is different from Example \ref{pos1}; since
$G(I)$ is \CM \ in  Example \ref{pos1}.
\end{example}

Before proving the theorem we give a proof of Theorem \ref{seccor}. We restate it 
for the convenience of the reader.
\begin{theorem}\label{as=2}
Let $(A,\m)$ be \CM \ local ring  of dimension $d = 3$. Let  $I$ be an $\m$-primary
ideal with $G(I)$ \CM \ and  $\red(I) \leq 2$. Then
\begin{center}
$f_2(I) \geq f_1(I) - f_0(I) + 1$.
\end{center}
\end{theorem}
\begin{proof}
As $G(I)$ is \CM \ we have that $\aF(I)  =  \red(I) -3 \leq -1$.  Let $J$ be a minimal reduction of $I$ such that $\red_J(I) = \red(I)$. 
As $G(I)$ is \CM \ the pair $(I,J)$
satisfy $V_{2}^{\infty}$. So the result follows from
Theorem \ref{sectheorem}.
\end{proof}

We now prove Theorem \ref{sectheorem}. The notation is  as in \ref{scrap}.
\begin{proof}
We use the complex $\D[n]$. By \ref{scrap2} and (\ref{hom}) we get $  $ 
\begin{equation}\label{r2}
-1 + 3\ell\left(\frac{I^n}{\m I^n}  \right) - 3\ell\left( \frac{I^{2n}}{\m I^{2n}} \right) +
 \ell\left( \frac{I^{2n}J^{[n]}}{\m I^{2n}J^{[n]} } \right) = \ell\left( H_2(\D[n]) \right), \ \  \forall n \gg 0.
\end{equation}
Since $\aF(I) < 0$ we have by (\ref{asymp}) $\red_{J^{[n]}}(I^n) = 2$ for all $n \gg 0$. So
$I^{2n}J^{[n]} = I^{3n}$ for all $n \gg 0$. Also for all $n \gg 0$ we have
$f_I(n) = \ell(I^n/\m I^n)$.
Setting these in (\ref{r2}) we get
$$-1 + 3f_I(n) - 3f_I(2n) + f_I(3n) \geq 0.$$
Since,
 $$ \quad f_I(n) = f_0\binom{n+2}{2} - f_1(n+1) + f_2,$$
an easy computation yields
\[
-1 + 3f_I(n) - 3f_I(2n) + f_I(3n) = -1 + f_0 - f_1 + f_2
\]
and the result follows.
\end{proof}

\begin{remark}
 Our second result Theorem \ref{sectheorem} has three hypothesis, namely 
\s \label{hypoth}
\begin{enumerate}[\rm (a)]
 \item 
$\grade(I) = l(I)$.
\item 
The pair $(I,J)$ satisfies $V_{2}^{\infty}$.
\item 
$\aF(I) < 0.$
\end{enumerate}
\end{remark}

 We show that if any of the hypothesis in \ref{hypoth} is not satisfied then the conclusion of Theorem \ref{sectheorem}
(i.e., $f_2(I) \geq f_1(I) -f_0(I) + 1 $) need not hold.

In the  first example only hypothesis \ref{hypoth} (a) is not satisfied, infact we have $l(I) - \grade(I) = 1$.
\begin{example} \label{ex1}
 Let $A = k[[X,Y,U,V]]/(XY, Y^3) = k[[x,y,u,v]]$ and $I = \m = (x,y,u,v)$. 
One can readily see that $\grade(\m) = 2$ while $l(\m) =  \dim A = 3$. Set $J = (x,u,v)$. Then  $I^3 = JI^2$.  So by 
\cite[3.2]{TrungProc} we get $\aF(I) <0$. 
The pair $(I,J)$ satisfy $V_{2}^{\infty}$
 by Proposition \ref{ex1J}. In fact in \ref{ex1J} we show   $\m^{2n} \cap J^{[n]} = J^{[n]} \m^n$ for all $n \geq 1$. 
 However 
 Hilbert series of $F(\m) = G(\m)$  (by COCOA) is 
\[
 \frac{1+ z -z^3}{(1-z)^3}
\]
So $f_0(I) = 1, f_1(I) = -2$ and $ f_2(I) = -3$. Thus $f_2(I) \ngeq f_1(I) -f_0(I) + 1$.
\end{example}
\begin{proposition}\label{ex1J}[With $A$, $\m$, $J$ as in \ref{ex1}]
\[
 \m^{2n} \cap J^{[n]} = J^{[n]} \m^n \quad \text{for all }  \ n \geq 1.
\]
\end{proposition}
\begin{proof}
 Fix $n \geq 1$. We  claim 
\begin{equation}
 \text{if}  \  x^nc \in \m^{2n}  \quad  \text{ then} \  c \in \m^n.  \tag{$\dagger$}
\end{equation}

Let us assume ($\dagger$)  and prove our result. Fix $n \geq 1$. Let $\xi \in \m^{2n} \cap J^{[n]}$. Write
$\xi = \alpha x^n + \beta u^n + \gamma v^n$. Set $(T,\q) = (A/(x^n), \m/(x^n))$. Note
\[
 T = \frac{k[[x,y,u,v]]}{(xy,y^3, x^n)} = \left(\frac{k[[x,y]]}{(xy,y^3, x^n)}\right)[[u,v]]
\]
Thus $u^*, v^*$ are $G(\q)$-regular. So 
\[
\ov{\beta}\ov{ u}^n + \ov{\gamma} \ov{v}^n \in  \ov{\m}^{2n} \cap (\ov{ u}^n,\ov{v}^n) = \ov{\m}^{n}(\ov{ u}^n,\ov{v}^n).
\]
It follows that $\ov{\beta}, \ov{\gamma} \in \ov{\m}^{n}$. Thus
\[
 \xi = \alpha x^n + \beta u^n + \gamma v^n + \theta x^n \quad \text{where} \ \beta,\gamma  \in \m^n \ \& \ \theta \in A.
\]
Notice $(\alpha + \theta )x^n \in \m^{2n}$. So by ($\dagger$) we get $(\alpha + \theta ) \in \m^n$. It follows that
$\xi \in \m^{n} J^{[n]}$.

We now
prove ($\dagger$).
Let $S =k[[x,u,v]]$ considered as a subring of $A$. 
Any element $a \in A$  can be written as
\begin{align} \label{AS}
 a &= \phi_{0}^{(a)}(x,u,v) + \phi_{1}^{(a)}(x,u,v)y + \phi_{2}^{(a)}(x,u,v)y^2  \\
\text{where} \ &\phi_{i}^{(a)}(x,u,v) \in S \ \text{for} 
\ i = 0,1,2.
\end{align}
So $A$, as a $S$-module, is generated by $1,y,y^2 $. Thus $\dim S = \dim A = 3$. It follows that  $S \cong k[[X,U,V]]$. 

Notice $A/(y) = S$  and the natural map $\pi \colon A \rt S$ is a splitting (as
$S$-modules) of the inclusion $\imath \colon S \rt A$. 
Set $L = Sy + Sy^2$. Then $A = S \oplus L$ as  a $S$-module. It follows that $\phi_{0}^{(a)}(x,u,v)$ in (\ref{AS}) is uniquely
determined by $a$. 

Let $\n$ be the unique maximal ideal of $S$. Then
\begin{equation} \label{AS2}
 \m^i = \n^i \oplus (\m^i \cap L) \quad \text{for all } \ i \geq 1.
\end{equation}

Set
 $$c = \phi_{0}^{(c)}(x,u,v) + \phi_{1}^{(c)}(x,u,v)y + \phi_{2}^{(c)}(x,u,v)y^2.$$ 
\textit{Notice}  $x^nc = x^n\phi_{0}^{(c)}(x,u,v)$. By uniqueness of $\phi_{0}^{*}$ we get
\[
 \phi_{0}^{(x^nc)}(x,u,v) = x^n\phi_{0}^{(c)}(x,u,v).
\]
Since  $x^nc \in \m^{2n}$, we get  by (\ref{AS2}) that $x^n\phi_{0}^{(c)}(x,u,v) \in \n^{2n}$. 
Clearly $x^*$ is $G_{\n}(S)$-regular. So $ \phi_{0}^{(c)}(x,u,v) \in \n^{n}$. By \ref{AS2} again we
get that $c \in \m^n$.

This proves $(\dagger)$. As stated earlier this finishes the proof of the proposition. 
 \end{proof}

In the  second example only hypothesis \ref{hypoth} (b) is not satisfied,
We adapt an example from \cite[6.2]{CJ3}. If $K$ is an ideal in $A$ let 
$\wt{K} = \cup_{n\geq 1} (K^{n+1} \colon K^n)$ be the Ratliff-Rush closure of $K$.
\begin{example}\label{ex2}
 Let $ A = \mathbb{Q}[[X,Y,Z]]$. Let $I = (X^4, X^3Y, XY^3, Y^4, Z)$. The ideal $J = (X^4, Y^4, Z )$ is a
minimal reduction of $I$, infact $I^3 = JI^2$. So by \cite[3.2]{TrungProc} we get $\aF(I) < 0$. Set
$B = \mathbb{Q}[[X,Y]]$ and $q = (X^4, X^3Y, XY^3, Y^4)$. One can show $\wt{q} \neq q$. However notice
$G(I) = G(q)[Z^*]$. So $Z^*$ is $G(I)$-regular. In particular $\wt{I} = I$. By
\cite[7.9]{Pu3} we get $\depth G(I^n) = 1$ for all $n \gg 0$. 
The Hilbert series of $F(I)$ is 
\[
 \frac{1 + 2z + 2z^2 - z^3}{(1-z)^3}.
\]
So $f_0(I) = 4$, $f_1(I) = 3$ and $f_2(I) = -1$. Thus $f_2(I) \ngeq f_1(I) -f_0(I) + 1$.
\end{example}

In the third example  hypothesis \ref{hypoth}(c) is not satisfied. Instead of \ref{hypoth}(b) the hypothesis
\begin{equation*}
 I^{2n}\cap (x_{1}^{n}, x_{2}^{n}) = I^n(x_{1}^{n}, x_{2}^{n}) \quad \text{for all} \ n \gg 1. \tag{$\dagger$}
\end{equation*}
is satisfied.  Recall $J = (x_1,x_2,x_3)$ is a minimal reduction of $I$. The hypothesis ($\dagger$)  is equivalent to
$\depth G(I^n) \geq 2$ for all $n \gg 0$  \cite[2.4]{Elasymp}. 
\begin{example}\label{ex3}
 Let $A = k[X,Y,Z]_{(X,Y,Z)}$ and $I = (X^3, XY^4Z, XY^5, Z^5, Y^7)$. 
Set $u = Z^5, V = 5X^3 + 3Y^7$  and $w=  X^3 -3XY^4Z  +2Z^5$.  Set $J = (u,v,w)$.  Using COCOA
we can check that $I^6 = JI^5$. So $J$ is a minimal reduction of $I$. 
The Hilbert Series of 
$G(I)$, $G(I/(u))$, $G(I/(u,v))$ is
\begin{align*}
 H(G(I),z) &= (1-z)H(G(I/(u)),z)  = (1-z)^2H(G(I/(u,v)),z), \\
           &= \frac{77 + 15z + 8z^2 + 2z^3 + 2z^4 + z^5}{(1-z)^3}, \   \text{and} \\
 H(G(I/(u,v,w)),z) &= 77 + 28z
\end{align*}
So $u^*, v^*$ is a $G(I)$-regular sequence. Note that $\depth G(I) = 2$. It follows that $\depth G(I^n) \geq 2$ for all
$n \geq 1$. So hypothesis ($\dagger$)  is satisfied.

 We prove that $\aF(I) \geq 0$.
Set $G = G(I)$.
\[
 a_{i}(I) = \max\{ n \mid H^i(G)_n \neq 0 \} \quad i = 0,1,2,3
\]
\s \label{compt} Note that $a_0(G) = a_1(G) = - \infty$. As $\red_J(I) = 5$ we get by \cite[3.2]{TrungProc}
that $a_3(G) \leq 2$. By \cite[2.1(a)]{Mar}  we have $a_2(G) < a_3(G) = 2$. 

By \cite[2.4]{Hoa} we have
\begin{equation*}
a_{i}(I^n)  \leq \left[ \frac{a_i(I)}{n} \right] \quad \text{for} \ i = 0,1,2.  \quad \& \quad
 a_{3}(I^n) = \left[ \frac{a_3(I)}{n} \right]
\end{equation*}
Notice by \ref{compt} we get 
\begin{align*}
 a_{i}(I^3) &\leq 0 \quad \text{for} \ i = 0,1,2. \\
a_3(I^3) &\leq -1 \quad \text{if }  \ a_3(I) < 0 \quad \text{and} \quad 
a_3(I^3) = 0  \quad \text{if }  \ a_3(I) \geq 0. 
\end{align*}
By \cite[3.2]{TrungProc} it follows that 
\begin{equation*}
 \red_{J^{[3]}}(I^3) = 2  \quad \text{if }  \ a_3(I) < 0   \quad \text{and} \quad 
\red_{J^{[3]}}(I^3) = 3  \quad \text{if }  \ a_3(I) \geq 0. 
\end{equation*}
However using COCOA we have verified that $I^9 \neq I^6J^{[3]}$. Thus
$\aF(I) = a_3(I) \geq 0$. The fiber coefficients are
$ f_0(I) = 17, f_1(I) = 34, f_2(I) = 17.$
So
\[
  f_{2}(I)- f_{1}(I) + f_{0}(I) - 1 = -1
\]
\end{example}
\section{Results when analytic spread is high }
In this section we extend our results Theorem \ref{result1} and
Theorem  \ref{sectheorem} 
to the case when $l(I) \geq 3 $ and  $l(I) \geq 4 $ respectively.
The main tool is the use of Rees-superficial sequences (see \ref{def}).
The utility of a Rees-superficial element in the study of fiber cone was 
first demonstrated in \cite{jpv}.

\s\label{def} An element $x \in I$ is said to be 
\emph{Rees superficial} if there exists  $r_0 \geq 1$ such that
$$(x) \cap I^r\m^s = xI^{r-1}\m^s \quad \text{for all} \ r \geq r_0 \ \& \ s \geq 0.$$ 
 
The following was proved  in \cite[2.8]{jpv} for $\m$-primary $I$ in a local ring $A$.
The same proof works in  general.
\begin{proposition}
 \label{diff} 
Let $(A,\m)$ be a local ring and $I$ an ideal.  Let $x \in I$
be a nonzero divisor in $R$ which is also Rees-superficial for  $I.$ 
Set $(B,\n) = (A/(x), \m/(x))$ and $K = I/(x)$. Then
$$ \ell\left(\frac{I^n}{\m I^n} \right) -  \ell\left(\frac{I^{n-1}}{\m I^{n-1}} \right) = \ell\left(\frac{K^n}{\n K^n} \right)  \quad \text{for all} \ n\gg 0.   $$
In particular $f_i(K) = f_i(I)$ for $i = 0,\ldots, l(I)-2$.
\end{proposition} 

The existence of a  Rees-superficial element which is also regular
follows from the    the following special case of a  Lemma (see \cite[1.2]{re}) due to Rees. 

\begin{lemma} \label{reeslemma}
Let $(A,\m)$ be  local and let $I$ be an ideal in $A$. Let $\mathcal{P}$ 
be a finite set of primes not containing $I\m$. Then
there exists $x \in I$ and 
 $r_0 \geq 1$ such that
\begin{enumerate}[\rm (1)]
\item
$x \notin \mathfrak{P}$ for all $\mathfrak{P} \in \mathcal{P}$ 
\item
$\displaystyle{(x) \cap I^r\m^s = xI^{r-1}\m^s \quad \text{for all} \ r \geq r_0 \ \& \ s \geq 0.}$
\end{enumerate}
\end{lemma} 

\begin{remark}\label{ReCsup}
If $x \in I$ is Rees-superficial and a non-zero divisor then 
it is easy to check that $x$
is $I$-superficial.
\end{remark}

 \s We say $x_1,\ldots,x_r \in I$ is a \emph{Rees-superficial} sequence
if $x_i$ is Rees superficial for the $A/(x_1,\ldots,x_{i-1})$-ideal $I/(x_1,\ldots,x_{i-1})$
for $i = 1,\ldots, r$.

\begin{remark}\label{rseq}
\textbf{a}. If $\grade(I) \geq r$ then using \ref{reeslemma} we can show  that there
exists a Rees-superficial sequence $x_1,\ldots,x_r$ in $I$
which is also a regular sequence. 

\textbf{b.} In this case we can further prove (by
using Proposition \ref{diff} repeatedly) that if $K = I/(x_1,\ldots,x_r)$ then
\begin{equation}\label{seq}
f_i(K) = f_i(I) \quad \text{for} \ i =  0,\ldots, l(I)-r -1.
\end{equation}
\end{remark}

We state our main results.
\begin{theorem}
\label{mth1}
Let $(A,\m)$ be   local and let $I$ be an ideal in $A$ with $l = l(I) \geq 2$ and
 $\grade(I) \geq  l(I)-2$. Assume either  $\reg^2(G(I)) \leq 1$ or $\red(I) \leq 1$. Then $f_1(I) \leq f_0(I) -1$. Furthermore equality holds if $\grade(I) =  l(I)$.
\end{theorem}
\begin{proof}
If $l(I) = 2$ then we are done by Theorem \ref{result1}. If $l \geq 3$ and 
as $\grade(I) \geq  l(I)-2$
   by \ref{rseq}.a we can choose a 
Rees-superficial sequence $x_1,\ldots,x_{l-2}$ which is also 
a regular sequence 
(and so an $I$-superficial sequence).
 Set $K = I/(x_1,\ldots,x_{l-2})$. Then $l(K) = 2$.

If $\reg^2(G(I)) \leq 1$
 then by \ref{notation} we get $\aF(K) + 2 = \reg^2(G(K)) \leq 1$. So
$\aF(K) < 0$. If $\red(I) \leq 1$ then $\red(K) \leq 1$. So by \cite[3.2]{TrungProc}
 we get $\aF(K) < 0$. Thus at any rate
$\aF(K) < 0$.

By Theorem \ref{result1} we get $f_1(K) \leq f_0(K) -1$. The result follows since by
 \ref{rseq}.b we have $f_i(I) = f_i(K)$ for $i = 0,1$.
\end{proof}

\begin{corollary}\label{NarC}
Let $(A,\m)$ be \CM \ local ring of dimension $d \geq 2$. If  $I$ is $\m$-primary and $e_{2}(I) = 0$ then
$f_1(I) = f_0(I) -1$.
\end{corollary}
\begin{proof}
 First assume that $d = 2$. By Narita's result \cite{Nar} we get that if $J$ is any reduction of $I$ then
$\red_{J^{[n]}}(I^n) = 1$ for all $n \gg 0$. Also as $\grade(I) = 2$ by Theorem \ref{result1}
we get $f_1(I) =  f_0(I) -1$.

When $d \geq 3$ we choose a Rees-superficial sequence $x_1,\ldots,x_{d-2}$ which is also an $A$-regular sequence (and so an $I$-superficial sequence). Set $K = 
I/(x_1,\ldots,x_{d-2})$. Note that $e_2(K) =  e_2(I) = 0$. Also $f_i(I) = f_i(K)$ for $i = 0,1$.
\end{proof}

Next we give an application of Theorem \ref{as=2}.

\begin{theorem}\label{mth2}
Let $(A,\m)$ be \CM \ local ring of dimension $d \geq 3$ and let $I$ be an $\m$-primary ideal. If $G(I)$
is \CM \ and $\red(I) \leq 2$ then
$f_2(I) \geq f_1(I) -f_0(I) +1$.
\end{theorem}
\begin{proof}
When $d = 3$ then the result follows from Theorem \ref{as=2}.
If $d \geq 4$ we choose a Rees-superficial sequence $x_1,\ldots,x_{d-3}$ which is also an $A$-regular sequence (and so an $I$-superficial sequence). Set $K = I/(x_1,\ldots,x_{d-3})$. Note that   $G(K)$ is 
\CM \ and reduction number of  $K$ is $\leq 2$.
By  Theorem \ref{as=2} we have $f_2(K) \geq f_1(K) -f_0(K) +1$. Also as 
$f_i(I) = f_i(K)$ for $i = 0,1,2$ (\ref{rseq}.b) we get the result.
\end{proof}

\section{appendix: Minimal reductions and Filter-regular elements}
In this section we prove that if $x^\circ  \in F(I)_1$ is $F(I)$ filter-regular then $(x^n)^\circ  \in F(I^n)_1$ is $F(I^n)$ filter-regular. This is used in  proof of Theorem \ref{veronese-th}. 

The relation between minimal reductions and Filter-regular sequences first appeared in the work of Trung
\cite{TrungProc}. We state one of  his results (\cite[3.8]{TrungProc}) in the form we need it. 

\begin{lemma}
 Let $I$ be an ideal with $s = l(I) > 0$. Then there exists a minimal reduction $J = (x_1,\ldots, x_s)$ of $I$ such that
$x_1^\circ, \ldots, x_s^\circ \in F(I)_1$ is a $F(I)$-filter regular sequence. 
\end{lemma}
  For definition of 
filter-regular sequence see \cite[18.3.7]{BH}. We however are only interested in a filter-regular element.

\s \label{filter-defn} Recall  that the following assertions  are equivalent
\begin{itemize}
 \item 
$x^\circ \in F(I)_1$ is $F(I)$-filter regular 
\item 
$(0 \colon_{F(I)} \ x^\circ)_n = 0$ for all $n \gg 0$.
\item 
 $x^\circ$ is $F(I)/H^{0}_{ F(I)_{+}}(F(I))$     regular.
\end{itemize}
For proof of the above equivalence see \cite[2.1]{TrungProc}. It is perhaps better to see Exercise \cite[18.3.8]{BSh}.
\begin{remark} \label{filtISreg}
By \ref{filter-defn} we get that if 
 $\depth F(I) > 0$ and $x^\circ$ is $F(I)$ filter-regular then $x^\circ$ is $F(I)$ regular.
\end{remark}

We give an ideal-theoretic criterion for an element $x^\circ \in F(I)_1$ to be $F(I)$-filter regular.

\begin{proposition}\label{filt-ideal-th}
Let $I$ be an ideal with $s = l(I) > 0$ and let $x \in I \setminus \m I$. \TFAE
\begin{enumerate}[\rm (i)]
 \item 
$x^\circ $ is $F(I)$ filter-regular.
\item 
$(\m I^{j+1} \colon x )\cap I^j = \m I^j$ for all $j \gg 0$.
\end{enumerate}
\end{proposition}
\begin{proof}
 (i) $\implies$ (ii)  We assume $(0 \colon x^\circ)_n = 0$ for all $n \geq c$. Clearly $\m I^j \subseteq (\m I^{j+1} \colon x )\cap I^j $ for all $j$.  If $a \in I^j \setminus \m I^j$ and $xa \in \m I^{j+1} $ then we have $x^\circ \bullet a^\circ = 0$. It follows that $j < c$. 

(ii) $\implies$ (i)  Coversely assume $(\m I^{j+1} \colon x )\cap I^j = \m I^j$ for all $j \geq c$. Say $a^\circ \in F(I)_j$ is non-zero and $x^\circ \bullet a^\circ = 0$. Then $a \in (\m I^{j+1} \colon x )\cap I^j$.  By hypothesis.
it follows that $j < c$. So $x^\circ $ is $F(I)$ filter-regular.
\end{proof}

\begin{corollary}\label{filt-ideal-cor}[with hypothesis as in \ref{filt-ideal-th}]
If $x^\circ  \in F(I)_1$ is $F(I)$ filter-regular then $(x^n)^\circ  \in F(I^n)_1$ is $F(I^n)$ filter-regular.
\end{corollary}
\begin{proof}
 Since $x^\circ$ is $F(I)$ filter-regular, by \ref{filt-ideal-th}, there exists $c > 0$ such that 
\[
 (\m I^{j+1} \colon x )\cap I^j = \m I^j \quad \text{ for all  }  \  j \geq  c.
\] 
So for $j \geq c $ we have $ (\m I^{j+n} \colon x^n )\cap I^j = \m I^j$. 
Therefore for  $j \geq c $ we obtain
\[
 (\m I^{n(j+1)} \colon x^n )\cap I^{nj} = (\m I^{(nj+n)} \colon x^n )\cap I^{nj} = \m I^{nj}
\]
Thus by \ref{filt-ideal-th} we get that $(x^n)^\circ$ is $F(I^n)$ filter-regular.
\end{proof}

 \providecommand{\bysame}{\leavevmode\hbox to3em{\hrulefill}\thinspace}
\providecommand{\MR}{\relax\ifhmode\unskip\space\fi MR }
\providecommand{\MRhref}[2]{%
  \href{http://www.ams.org/mathscinet-getitem?mr=#1}{#2}
}
\providecommand{\href}[2]{#2}

%\bibliographystyle{amsplain}
%\bibliography{cref}
\end{document}